\documentclass[12pt,notitlepage,twoside,a4paper]{amsart}
 \usepackage{amsfonts} 
\usepackage[latin1]{inputenc}
\usepackage[cyr]{aeguill}
\usepackage{xspace}
\usepackage[polutonikogreek,english]{babel}

\usepackage{amsmath,amssymb,enumerate}

\usepackage{epsfig,fancyhdr,color}

\usepackage{amssymb}
\usepackage{amsmath,amsthm}
\usepackage{latexsym}
\usepackage{amscd}
\usepackage{psfrag}
\usepackage{graphicx}
\usepackage[latin1]{inputenc}
\usepackage[all]{xy}
\usepackage[mathcal]{eucal}

\definecolor{NoteColor}{rgb}{1,0,0}

\renewcommand{\textsc}{\textcolor{red}}

\newtheorem*{theorem 1}{\rm\bf Proposition 1}
\newtheorem*{theorem 2}{\rm\bf Proposition 2}

\theoremstyle{definition}

\theoremstyle{remark}

\def\interieur#1{\mathord{\mathop{\kern 0pt #1}\limits^\circ}}

\title[Sullivan: Mathematics and physics]{Dennis Sullivan: Mathematics and physics; manifold and space}

\author{Athanase Papadopoulos}
 
 \date{\today}

\begin{document}
  
  \maketitle
  
  \begin{abstract}
  This article is a personal overview of the work of Dennis Sullivan who was awarded the 2022 Abel prize. It was commissioned by  the Bulletin of the (Indian) Mathematics Consortium, and it will appear there.
  \medskip
  
  \noindent
  Keywords: Abel Prize, Dennis Sullivan, topology, manifold, space, hydrodynamics.

    \medskip
  
  \noindent  AMS Classification:  01A70
  \end{abstract}
     
  We do not have a Nobel prize in mathematics, but we have two distinctions which are at least comparable to it, the Fields medal and the the Abel prize. Unlike the Fields medal, the Abel prize  has no age limitation; it is usually awarded for long-term achievements. 
  
   The history of the Abel prize is interesting. Sophus Lie, who, like Niels Henrik Abel, was Norwegian, proposed the establishment of this prize already in 1899, after he learned that the Nobel Prizes (which were to be awarded for the first time in 1901 by the Swedish and the Norwegian Academies) will not include mathematics. The first Abel prize was planned to be part of the events celebrating the 100th anniversary of  Abel's birth (August 5, 1802). The project was then delayed, and eventually aborted for political reasons, in particular the dissolution of the union between the kingdoms of Norway and Sweden (1905). The idea of this prize arose again at the beginning of the second millenium, and in August 2001, the Norwegian government announced that the prize would be awarded starting in 2002, again in relation with Abel, this time for his two-hundredth anniversary. The prize was slightly delayed, and it was attributed for the first time in 2003, to Jean-Pierre Serre, who was awarded the Fields medal about 49 years before (at age 27). It  was given every year since then, and the names of the recipients include Michael Atiyah and Isadore Singer (joint prize), S. R. Srinivasa Varadhan, Misha Gromov, John Milnor, Yakov Sinai, John Nash and Louis Nirenberg (joint), and Andrew Wiles, for their outstanding contributions in various fields of mathematics. The name of Dennis Sullivan was like a gap in this list, and it was filled in March 2022. 
     
  It is impossible to do justice in a few pages to the titanic work that Sullivan accomplished in a period of 60 years, to the exceptional influence he had on shaping 20th-21th century mathematics, and to all the good he has done for the mathematical  community. In this report, after a brief \emph{Vita}, I have chosen to start by reviewing  the three lectures that Sullivan gave at three International Congresses of Mathematicians (Nice 1970, Vancouver 1974 and Berkeley 1986). These lectures reflect some of his subjects of interest during the first part of his career. Then, I will present a list of the incredible amount of topics on which he worked,  with special emphasis on dynamics and  fluid mechanics.  From time to time I will quote some emails I have exchanged with him over the years. Going through this correspondence, I see in hindsight that he was always answering my questions, mathematical or not, and his remarks were always compelling. I will mention in particular some thoughts he shared recently with me and with several of his colleagues and friends on his view on what is important in mathematics.
  
 \medskip
  
  Born in 1941, in Port Huron (Michigan), Sullivan was brought up in Houston (Texas). He entered Rice University, first as a student in chemical engineering. It was the discovery of topology that made him change his mind and shift to mathematics. In an interview with 
Shubashree Desikan which appeared in \emph{The Hindu} on March 24, 2022, Sullivan speaking informally about this episode says: ``At Rice University, all the science students, electrical engineers and all the others, took math, physics and chemistry. In the second year, when we did complex variables, one day, the professor drew a picture of a kidney-shaped swimming pool, and a round swimming pool. And he said, you could deform this kidney-shaped swimming pool into the round one. At each point, the distortion is by scaling. A little triangle at this point goes to a similar triangle at the other point.  We had a formula for the mapping, because we were taking calculus, and we had a notation for discussing it. This was like a geometric picture. This mapping was essentially unique. The nature of this statement was totally different from any math statement I've ever seen before. It was, like, general, deep, and wow! And true! So then, within a few weeks, I changed my major to math."

  Sullivan obtained his PhD at Princeton University in 1966, with a dissertation titled \emph{Triangulating homotopy equivalences}, with William Browder as advisor. In the few preceding years, algebraic topology had been invigorated with the introduction of new techniques and with a series of outstanding results on the classification of manifolds. Let me say a few words on this period.
  
  The notions of fiber space and fiber bundle became central  around the year 1950. Soon after, Serre  introduced  in his thesis (1951) the crucial idea of using spectral sequences to study the homology of fiber spaces.  In 1952, Michail M. Postnikov obtained his famous  reconstruction of  the cohomology of a space from its homotopy invariants.
 René Thom, in 1954, obtained a classification of manifolds up to cobordism. In 1956, Milnor proved the existence of exotic differentiable structures on the 7-dimensional sphere.  The higher-dimensional analogue of the Jordan curve theorem (the so-called Jordan--Schoenflies theorem) was obtained in 1959-1960 (works of Barry Mazur, Marston Morse and Morton Brown).
 In 1960,  Steven Smale proved the Poincaré conjecture (PL and smooth categories) for dimensions $\geq 7$ (and then $\geq 5$). In 1961, Christopher Zeeman obtained a proof of the same conjecture for dimension 5, and in 1962 Stallings gave a proof for dimension 6.
  In 1961, Milnor disproved the so-called \emph{Hauptvermutung der kombinatorischen
Topologie} (``main conjecture of combinatorial
topology"), a conjecture formulated in 1908 by Ernst Steinitz and Heinrich Franz Friedrich Tietze, asking whether  any two triangulations of homeomorphic spaces are  isomorphic after subdivision.   
 In 1965 Sergei P. Novikov proved that the Pontryagin classes with rational coefficients of vector bundles are topological invariants. 
 
Sullivan started working on topology amid all this booming activity.
 In his thesis, he gave an obstruction to deforming a homotopy equivalence of piecewise linear manifolds to a PL homeomorphism, obtaining further cases where the Hauptvermutung is true. This result is called ``the characteristic variety theorem." It was used  to provide numerical
invariants that classified the combinatorial manifolds in a homotopic type. 
  He summarized his results on this topic in his short paper 
 \emph{On the Hauptvermutung for manifolds} \cite{S1967} which appeared in 1967 and for which he was awarded the Oswald Veblen prize of the AMS.

After his PhD, Sullivan worked successively  at the  University of Warwick (1966-1967), the University of California at Berkeley (1967-1969) and MIT (1969-1973). In the year 1973-74, he was invited as a visiting professor at the Université de Paris-Sud, Orsay. The next year, he became a permanent member of the Institut des Hautes \'Etudes Scientifiques at Bures-sur-Yvette. The institute is a couple of kilometers from the Orsay campus. 
     Remembering that period, I find in an email from Sullivan (2015): ``Grothendieck left IH\'ES around 1970. Quillen visited IH\'ES from MIT during the year 1972-1973. I visited IH\'ES  and Orsay from MIT during the year 1973-1974. It was a splendid place to do Math. IH\'ES  offered Grothendieck's vacated position to Quillen who declined.  IH\'ES offered it to me and I grabbed it."      
     
      In a 2019 email, Sullivan writes: ``My first math hero, as a grad student in Princeton, was René Thom. When IH\'ES offered me a professorship, 
I was honored to accept it and to become Thom's `colleague'." He adds: ``My second math hero starting in December 71  was the  Mozart-like  figure  Bill Thurston.
During the next decade I was in France and Bill was in the US  but we had a quite fruitful and intense  interaction." I shall say more about Sullivan's relation with Thurston below.

 In 1981, Sullivan was appointed to the Einstein Chair at the Graduate Center of the City University of New York.  He kept his position at IH\'ES on a part-time basis and started spending half the year in France and the other half in the US, until 1996, where he took up a professorship at the State University of New York at Stony Brook, again on a part-time basis, keeping his Einstein Chair at CUNY, where he conducts, since his appointment there, a weekly seminar on geometry in the broadest sense whose sessions are known to last for several hours (sometimes all day). For me, and I assume it is the same for many others who know this place,  the Graduate Center of CUNY is the central point where the heart of New York is beating.

To talk about Sullivan's work, I start with the lecture that he gave at the 1970 ICM, which took place in Nice.   The invitation to that Congress came just after Sullivan proved (concurrently with and independently of D. Quillen) the Adams conjecture, which concerns the homotopy theory of 
sphere bundles associated with vector bundles. The conjecture was considered as one of the most important conjectures in topology. Whereas this led Quillen to develop algebraic K-theory, Sullivan's approach was based on  the ``arithmetization'' of geometry and topology, more precisely on the introduction of Galois theory in the  geometry of manifolds in the form of concepts like localisation, rationalization, profinite and $p$-adic homotopy theory. Sullivan's proof of Adam's conjecture, which he obtained in 1967,   is based on the construction of a functor from abstract algebraic varieties into profinite homotopy theory. This approach led him to the study of the absolute Galois group itself through its actions on new geometric objects, and at the same time gave rise to strong relations between number theory and homotopy theory. In fact, Sullivan studied the action of the absolute Galois group in homotopy theory via classifying spaces and
Postnikov towers. The ICM talk is titled \emph{Galois symmetry in manifold theory at the primes} \cite{ICM70}. 
A more detailed version of this work appeared in the paper \cite{S1974}, in which Sullivan describes his (still open) unrequited \emph{Jugendtraum} (childhood dream).  Let me quote from the introduction to that paper some sentences which are characteristic of Sullivan's personal thinking, based on rich and inspiring analogies:
``We  are studying the structure of 
homotopy types to deepen our understanding of more complicated or richer
mathematical objects such as manifolds or algebraic varieties.
The relationship between these two types of objects is I think rather
strikingly analogous to the relationship in biology between the genetic structure
of living substances and the visible structure of completed organisms or
individuals.
The specifications of the genetic structure of an organism and of the
homotopy structure of a manifold have similar texture; they are both discrete,
combinatorial, rigid, interlocking and sequential."

Sullivan's second ICM talk (Vancouver 1974) is  titled : ``Inside and outside manifolds"  \cite{ICM74}.
The title is characteristic of Sullivan's  taste for literary turns of phrase. The paper has two parts, which express his two approaches at the time for the exploration of manifolds: topology and dynamical systems. The first part, on ``outside" manifolds, is an exposition of the complete classification theory of simply-connected manifolds of dimension at least 5. The main detailed reference for this part is a set of notes on geometric topology  based on lectures that Sullivan gave at MIT in 1970, carrying the subtitle ``Localization, periodicity, and Galois symmetry"  \cite{MIT-notes}. These notes were widely circulated and they had a major influence on algebraic and geometric topology not only in the West but also in the Soviet Union where they were translated into Russian just after Sullivan's ICM 1974 talk.
In the second part of the lecture, ``Inside manifolds", Sullivan presents a series of qualitative results on dynamics (in particular Smale's Axiom A, structural stability, and the asymptotic properties of leaves of foliations) within individual manifolds. This second part, as Sullivan writes, ``focuses attention on the classical goals and problems of \emph{analysis situs}."

A new version of the 1970 MIT notes was edited in 2005 by Andrew Ranicki  \cite{MIT-notes}. It concludes with a 10 pages postscript by Sullivan in which he recounts, in his personal literary style, the genesis and later developments of these notes, as well as episodes from his own mathematical education with the mathematicians who have inspired him. 
René Thom, who, as we recalled, was Sullivan's first hero, is abundantly quoted. The postscript contains several open problems and conjectures. At several points, the tone is philosophical. In his book \emph{Surgery on compact manifolds}, C. T. C. Wall writes about these notes: ``It is difficult to summarise Sullivan's work so briefly: 
the full philosophical exposition should be read."  The postscript is punctuated with details on Sullivan's family life.  After recalling his work in the 1970s, Sullivan writes in these notes:  ``About this time dynamical systems, hyperbolic geometry, Kleinian
groups, and quasiconformal analysis which concerned more the geometry of the manifold than its algebra began to distract me (see ICM report 1974), and some of the work mentioned above was left incomplete and unpublished."
This mention of quasiconformal mappings brings us to his 1986 ICM talk titled \emph{Quasiconformal homeomorphisms in dynamics, topology, and geometry}.

Sullivan's third ICM talk (Berkeley 1986) \cite{ICM86}  is a survey of results and conjectures, mainly due to him, on   quasiconformal homeomorphisms and their use in four different contexts. Here, a homeomorphism $\varphi:X\to Y$ between metric spaces $(X,\vert \  \vert)$ and $(Y,\vert \  \vert)$ is said to be  ($K$-)quasiconformal if for some $K>0$ 
$$\limsup_{r\to 0}\frac{\sup\vert \varphi(x)-\varphi(y)\vert \hbox{ where } \vert x-y\vert =r \hbox{ and $x$ is fixed }}{\inf\vert \varphi(x)-\varphi(y)\vert \hbox{ where } \vert x-y\vert =r \hbox{ and $x$ is fixed }} < K.
$$

The four contexts in which quasiconformal mappings are studied in this paper are the following:

\medskip

 \noindent (1) Feigenbaum's numerical discoveries in 1-dimensional dynamics. Here,  quasiconformal homeomorphisms are used to define a distance between real analytic dynamical systems, by first complexifying them. This distance  is contracting under the Feigenbaum renormalization operator.

 \medskip
 \noindent (2) The theory of quasiconformal manifolds, including developments of de Rham cohomology, Atiyah--Singer index theory and Yang--Mills connections on these manifolds. This is joint work between Sullivan and Donaldson. Combined with other works of Donaldson and Freedman, this theory provides a complete picture of the structure of quasiconformal manifolds, namely, each topological manifold has an essentially unique quasiconformal structure,  
 except in dimension 4 where, by results of Freedman and Donaldson--Sullivan,  the statement is false. 
 
 \medskip
 \noindent (3) The quasiconformal deformation theory of expanding systems $f:U_1\to U$ where $U_1$ is a domain of the Riemann sphere and $f$ a $d$-sheeted ($d>1$) onto covering.  Sullivan writes: ``The analytic classification of expanding systems of a given topological dynamics type on the invariant set is a kind of Teichm\"uller theory." Here, an infinite-dimensional Teichm\"uller space of expanding analytical dynamical systems near their fractal invariant sets is embedded in the Hausdorff measure theories possible for the transformation on the fractal, and the Hausdorff measure theories of fractals are embedded in the theory of Gibbs states.

 \medskip
 \noindent (4)   The quasiconformal theory of the geodesic flow of  negatively curved manifolds via the action of the fundamental group on the sphere at infinity of the universal cover. Constant curvature is given a characterization among variable negative curvature in terms of a uniform quasiconformality property of the geodesic flow.
 
\medskip
  
  These three ICM talks give us an idea of part of Sullivan's mathematical interests until the year 1986.   I will mention some other results of him, but before that I would like to list all the fields and the subjects on which he worked. 
 
 We already mentioned  his work on the topology of manifolds. His results on this topic concern, among others, the classification of manifolds in several categories: smooth, PL, topological, Lipschitz, bi-Lipschitz and quasiconformal.  Closely related is his work on algebraic topology: homology and  homotopy  theories, including rational homotopy theory, intersection homology, differential cohomology, homotopical algebra, characteristic classes, K-theory,  flat bundles, minimal models, Galois symmetry and string topology. Then come the more geometrically specialized topics: foliations, minimal surfaces,  geometric structures and geometry of manifolds: affine manifolds, complex manifolds, K\"ahler manifolds, hyperbolic geometry, Kleinian groups, universal Teichmüller theory, laminations, the soleniod and circle packings. Under the heading ``dynamics", I mention 1-dimensional dynamics, renormalization, universality, holomorphic dynamics, potential theory, optimal control, measurable dynamics, and the related topics of ergodic theory, probability, Brownian motion, Diophantine approximation,   chaos and fractals (in particular, self-similar structures arising in the theory of Kleinian groups, in KAM theory and in iteration theory of holomorphic dynamics).
Finally, let me also mention in a nutshell:  $C^*$ algebras, harmonic analysis,
 singularity theory, field theories, sigma models, and logic.

 I do not know of any mathematician who has such a broad spectrum of interests, except Leonhard Euler. 
  
  \medskip

  I would like to say a few more words on Sullivan's French period.
  
  As a visitor at Orsay, Sullivan gave a course on the theory of rational homotopy of differential forms that he had newly developed. During the same period, he came with a new interpretation of the Godbillon--Vey invariant, a subject that was keeping busy geometers working on foliation theory which was a hot topic at that time at Orsay. The interpretation was completely new, using the notion of currents. Sullivan included this result in his paper \emph{Cycles for the dynamical study of foliated manifolds and complex manifolds} \cite{Sullivan76}, which appeared in 1976 and which contains a wealth of other results. It was also during the same year that he introduced the topologists at Orsay to Thurston's theory of surface homeomorphisms; this was the major impetus to the seminar on Thurston's work which took place at Orsay in 1976-77 and which gave rise to the famous book ``Travaux de Thurston sur les surfaces."
  
  Talking about Orsay and Bures, and on a more personal level, let me recall an episode, in 1983. I was a PhD student at Orsay. François Laudenbach, who was my thesis advisor, had asked me a question related to pseudo-Anosov homeomorphisms. I came to see him one day, with a four pages manuscript, containing a proof of a theorem which answers the question. He told me: ``I will show it to Sullivan." After he talked with Sullivan, he decided that this was my PhD dissertation. (For the record, I had some trouble later on, convincing people at the registrar's office that these 4 pages constitute a thesis dissertation.)

  I already quoted Sullivan saying that he had two heroes, Thom and Thurston, and I would like to say a few more words on his relation with the latter. 

Sullivan was always a major promoter of Thurston's ideas, and he was probably
the person who best understood, from the first years, their originality and their broadness. The main themes  discussed at Sullivan's seminar at IH\'ES  and at  his New York Einstein Chair seminar included 1-dimensional dynamics and the so-called kneading theory established by Milnor and Thurston. The other topics include Kleinian groups (discrete isometry groups of hyperbolic 3-space) and holomorphic dynamics, two topics which eventually became a single topic, after Sullivan established a complete dictionary between them. Here, again,
Thurston's ideas were often at the forefront, and Sullivan spent years trying to understand and to explain them.
 He  was the first to learn from Thurston his result on the characterization of  postcritically finite rational maps of the sphere, that is, rational maps whose
forward orbits of critical points are eventually periodic. The proof of this theorem, like
the proofs of several of Thurston's major theorems, uses a fixed point argument for an action on a Teichmüller space. The rational map in the theorem is obtained
through an iterative process as a fixed point of that map on Teichm\"uller space. This became a major element in  Sullivan's analogies between the iteration theory of rational maps and the theory of Kleinian groups.  Thurston's theorem, together with Sullivan's dictionary between the theory of discrete subgroups of $\mathrm{PSL}(2,\mathbb{C})$ and complex analytic iteration, constitute now the two most
fundamental results in the theory of iterations of rational maps. Regarding this theory, in 1985, Sullivan published a paper in which he gave the proof of a longstanding
question formulated by Fatou and Julia in the 1920s and which  became known
as the Sullivan's \emph{No-wandering-domain Theorem}. This theorem says that every component of the Fatou set
of a rational function is eventually periodic. 
A fundamental tool that was introduced by
Sullivan in his proof is that of quasiconformal mappings, the main topic of his 1986 ICM lecture and one of the major concepts used in
classical Teichmüller theory.

 Let me also mention that in the realm of conformal geometry, Thurston introduced the subject of discrete
conformal mappings, and in particular the idea of discrete Riemann mappings. In 1987,
Sullivan, together with Burton Rodin, proved an important conjecture of Thurston on
approximating the Riemann mapping using circle packings. This result, the authors write, is in the setting of Thurston's ``provocative, constructive, geometric approach  to the Riemann mapping theorem," see  \cite{RS1987}.

 There is one field which I did not stress upon yet; this is physics.
  
Sullivan is also a physicist, in the tradition of Euler and Riemann, who were also physicists. In fact, Riemann considered himself more a physicists than a mathematician.  He was thoroughly involved in  gravitation, electricity, magnetism and electrodynamics and he adopted a completely physical approach to  the theory of functions of a complex variable.  In a note contained in his \emph{Collected Works}, he writes: ``My main work consists in a new formulation of the known natural laws---expressing them
in terms of other fundamental ideas---so as to make possible the use of experimental data on the exchanges between heat, light, magnetism, and electricity."    During his stay at IH\'ES, Sullivan was discussing with physicists, Oscar Lanford, Henri Epstein, David Ruelle, etc.
 In the Postface to the 2005 edition of his 1970 MIT notes, he writes:  ``At a physics lecture at IHES in 1991 I learned the astonishing (to
me) fact that the fundamental equations of hydrodynamics in three
dimensions were not known to have the appropriate solutions."

This brings us to hydrodynamics, one of Sullivan's favorite topics, and in particular, the  Navier--Stokes equation. This is the main partial differential equation which describes the flow of an incompressible fluid. It is good to remember here that a special form of this equation was discovered by  Euler back in the 18th century (Euler's equation is special in that it does not take viscosity into account).  In the interview that appeared in \emph{The Hindu} which we already quoted, Sullivan, talking about the Navier--Stokes equation, declares that
his fascination for this subject goes back to his youth spent in Texas. He says, again with his informal tone: 
``If you're in Texas, as a student of chemical engineering, there's the petrochemical industry, the oil industry, and organic chemistry and plastics, all around Houston. If you are good in science, and you work on that and become an engineer you can get a good job and have a nice work at a research center. So it's a good thing to do. In fact, during the summers, I had jobs at various such places. Once I had to study the computer methods that they were using to do what's called secondary recovery. You know, when they find oil, because of the pressure, if they drill a hole, the pressure makes it shoot up, right. But after they drill for 20 years, the pressure goes down. What they do then is go to another part of the field, and they drill and they pump in water to create pressure that will push the oil back to their wells, and for this they have to solve the linearised version of the Navier--Stokes equation. I didn't know that name, then but it's a linearised version of the Navier--Stokes equation. While at the summer job where I was studying the possible computer programs I had a certain question there. That was around 1960. [\ldots] So in a sense, I was aware that there's this huge industry related to fluid flow through porous media. It was astonishing to me to find out as I found out in the 1990s, that that equations in three dimensions, the beautiful equations, are not solved."

In a recent email (2021), Sullivan writes on the same subject: ``Back to 1991-92:  After a talk at  IHES  by  computational physicists using Fourier modes to  numerically analyze  the 3D incompressible  Navier--Stokes Equation one learned about the lack of a mathematical answer to the question of long time existence of solutions in certain classes. That  this beautiful and  widely used equation that designed airplanes helped to more efficiently transport  oil and  to understand how to stabilize aneurisms in or near the brain or heart etc. held such mysteries was astonishing.   For example, why  was this 3D problem so hard? The same Navier--Stokes  problem in 2D could be treated by exactly the same  Calderon--Zygmund analysis tools that treated the Beltrami equation, so important in the Ahlfors--Bers treatment of Teichm\"uller theory.  In addition half of the non-trivial theory of quasiconformal mappings in 2D extended to higher dimensions as in Mostow rigidity for  3D."

Let me also recall that the Navier--Stokes equations are among the Clay Mathematics Institute Millennium Prize problems.  The statement of the Clay problem is the following:
``Waves follow our boat as we meander across the lake, and turbulent air currents follow our flight in a modern jet. Mathematicians and physicists believe that an explanation for and the prediction of both the breeze and the turbulence can be found through an understanding of solutions to the Navier--Stokes equations. Although these equations were written down in the 19th Century, our understanding of them remains minimal. The challenge is to make substantial progress toward a mathematical theory which will unlock the secrets hidden in the Navier--Stokes equations."

Talking about physics, let me also mention that Sullivan's work on 1-dimensional dynamics, period doubling and chaos was motivated by physical experimentation. 
 In the Postface to his MIT notes, he writes: ``For example, one new project began in the 1980s in fractal geometry---the Feigenbaum universal constant associated to period doubling---presents a new kind of epistemological problem. Numerical calculations showed that a certain mathematical statement of
geometric rigidity in dynamics was almost certainly true, but the
available mathematical technique did not seem adequate. Thus assuming the result was true there had to be new ideas in mathematics
to prove it. The project consumed the years up to the birth of my
second son Thomas in July `88."

 Recently (2021), Sullivan, emailed emailed me together with a few friends and colleagues some basic and fundamental questions and thoughts on ``What is a manifold", describing this as his mathematical lifelong quest. I read in one of these emails: ``3D fluids and turbulence (as in the real world)  is an area that has been puzzling me since the 90-91 because  it seems its natural setting or comfort zone is not really identified yet. I am thinking, focusing really on  \emph{What is a manifold}". In another mail, I read:   ``Riemann introduced the notion of manifold, Gelfand, then Grothendieck
and then Connes recast the geometrical aspect in terms of the algebra and the quality of functions. The two viewpoints are in duality." 

A few months ago, he sent me to read a 1992 paper by  Hikosaburo Komatsu on hyperfunctions and microfunctions, with the comments:  ``Here is a marvelous example of a transfer of information explaining the history of ideas. [\ldots] 
I liked so much how the  paper  first explained the ease of the 1-dimensinal case (in terms of what we know from before 1900), and then the modern story with the competitions between wave front set perspective of  one school and the sheaf/algebra perspective of another. [\ldots] I  think this  sheaf of hyperfunctions is a key point for me in a three decade  hiatus to define smooth structures in real dimension one that  will accommodate the vastly successful  theory of dynamics in one dimension
which provides the tool and the language to  be able to `explain' some of the missing parts of  the numerical universalities discovered numerically by physicists in the 70's . The non-missing part achieved by the end of the 80's only worked for integer critical points using Teichm\"uller theory, even though the numerical theorems were known for all critical exponents."   In another email, again, on the same topic: ``The thread  with which I am preoccupied is the discussion of space where physical reality seems to take place which, by the way, can be reformulated  with charts in terms of the notion of function or with algebra in terms of ideals in an algebra of functions".

We are back to the most fundamental questions: What is a manifold? What is a function? What is space?

\end{document}